\documentclass[12pt]{amsart}

\usepackage{amsmath}
\usepackage{amssymb}
\usepackage{latexsym}
\newtheorem{theorem}{Theorem}[section]

\newtheorem{lemma}[theorem]{Lemma}
\newtheorem{prop}[theorem]{Proposition}
\newtheorem{defn}[theorem]{Definition}

\newenvironment{proof*}{\vskip 2mm\noindent {}}{\hfill $\Box$ \vskip 2mm}
\numberwithin{equation}{section}
\newcommand{\C}{{\mathbb{C}}}

\newcommand{\R}{{\mathbb{R}}}

\newcommand{\D}{{\mathbb{D}}}

\newcommand{\eps}{\varepsilon}

\begin{document}

\title{Convergence and multiplicities for the Lempert function}



\author{Pascal J. Thomas, Nguyen Van Trao}
\maketitle

\begin{abstract} Given a domain $\Omega \subset \mathbb C^n$, 
the Lempert function is
a functional on the space $Hol (\D,\Omega)$ of
analytic disks with values in $\Omega$, depending on a set of poles
in $\Omega$. We generalize its definition to the case where poles
have multiplicities given by local indicators (in the sense of Rashkovskii) 
to obtain a function which still dominates the corresponding Green
function, behaves relatively well under limits, and is monotonic with
respect to the local indicators. In particular, this is an improvement over
the previous generalization used by the same authors to find an
example of a set of poles in the bidisk so that the (usual) Green and
Lempert functions differ. 
\end{abstract}

\section{Introduction}
\label{intro}
We assume throughout that $\Omega $ is a bounded domain in $\C^n$. Let $\D$ stand for the unit disk in
$\C$. 
The classical Lempert function with pole at $a \in \Omega $
\cite{Lempert} is defined by
$$
\ell_a (z):=\inf \big\{ \log|\zeta|:\exists \varphi\in Hol (\mathbb
D,\Omega),
\varphi(0)=z,  \varphi(\zeta)=a\big\}.
$$

Given a finite number of points $a_j \in \Omega$, $j=1,...,N$,
Coman \cite{Coman} extended this to:
\begin{multline}
\label{coman}
\ell (z):=\ell_{\{a_1,\dots,a_N\}} (z):=
\inf \big\{ \sum^N_{j=1} \log|\zeta_j|:
\\
\exists \varphi\in Hol (\mathbb D,\Omega) :  \varphi(0)=z,
\varphi(\zeta_j)=a_j, j=1,...,N \big\} .
\end{multline}

The \emph{Green function} for the same poles is
\begin{multline*}
g := \sup  \left\lbrace u \in PSH(\Omega, \mathbb R_-) :  u(z) \le
\log |z-a_j|+C_j, \right.
\\
\left.
\mbox{ for } z \mbox{ in a neighborhood of } a_j, j=1,...,N \right\rbrace ,
\end{multline*}
where $PSH(\Omega, \mathbb R_-)$ stands for the set of all negative
plurisubharmonic functions in $\Omega$.
The inequality $g(z)\le \ell(z)$ always holds, and it is known that it
can be strict
\cite{CarlWieg}, \cite{TraoTh}, \cite{NikoZwo}. If $\ell$ ever turns out
to be plurisubharmonic itself, then $\ell$ must be equal to $g$
\cite{Coman}.

There are natural extensions of the definition of the Green function. 
In one dimension, considering a finite number
of poles in the same location $a$, say $m$ poles, has a natural
interpretation in terms
of multiplicities:
the point mass in the Riesz measure of the Green function is multiplied
by $m$. Locally, the Green function behaves like $\log |f|$, where
$f$ is a holomorphic function vanishing at $a$ with multiplicity $m$.

Lelong and Rashkovskii \cite{lelongrash}, \cite{Rash} defined a
generalized Green
function. The function $\log |z|$ was replaced by 
``local indicators",
 i.e. circled plurisubharmonic functions $\Psi$ whose Monge Amp\`ere
measure $(dd^c \Psi)^n$ is concentrated at the origin, such that whenever $\log |w_j|=c\log |z_j|$
for all $j \in \{1,\dots,n\}$, then
$\Psi(w)=c\Psi(z)$. This
has the
advantage of allowing the consideration of non-isotropic singularities such as
$\max(2\log |z_1|, \log |z_2|)$, but the ``circled" condition privileges
certain coordinate axes, so that the class isn't invariant under
linear changes of variables. We will have to remove this restriction to
obtain a class large enough to describe some natural limits.

In several complex variables, we would like to know which notion of
multiplicity can arise when we take limits of ordinary Green
(or Lempert) functions with several poles tending to the same point.
This idea was put to use in \cite{TraoTh} to exhibit an example where a
Lempert function with four poles is different
from the corresponding Green function. The definition of a generalized Lempert function
chosen in \cite{TraoTh} 
had some drawbacks --- essentially, it was not monotonic with respect to
its system of poles (in an appropriate sense)
 \cite[Proposition 4.3]{TraoTh}
 and did not pass to the limit in some very simple situations
\cite[Theorem 6.3]{TTppt}. We recall that monotonicity holds
when no multiplicities
are present,
see \cite{WikstromAMS} and \cite[Proposition 3.1]{TraoTh} for the convex
case, and the more recent \cite{NiPfl} for arbitrary domains
and weighted Lempert functions, or more generally when a subset 
of the original set of poles is considered with the same generalized local indicators.

In section \ref{definition}, we successively
define a class of indicators, a subclass which is useful to produce ``monomial"
examples, a notion of multiplicity for values attained
by an analytic disk, and a generalization of Coman's Lempert function
to systems of poles with generalized local indicators, different from 
\cite{TraoTh}. In section \ref{main}, we state our two main results:
monotonicity, and convergence under certain restrictive (but, we
hope, natural) conditions. Further sections are devoted to the proofs
of those results.

Finally, in Section \ref{compprev} we summarize the differences between
our new definition and that given in \cite{TraoTh}.

The first named author would like to thank Nikolai Nikolov for stimulating discussions
on this topic, and his colleague Anne Bauval from Toulouse for showing him
a nice purely combinatorial proof of Lemma \ref{order}. Special thanks are due
to the referee for his very thorough reading of our paper. 

\section{Definitions}
\label{definition}
\begin{defn} 
\cite{lelongrash}
\label{LRindic}
Let $\Psi \in PSH(\D^n)$. We call $\Psi$ a \emph{local
indicator}
and write $\Psi \in \mathcal I_0$ if
\begin{enumerate}
\item
$\Psi$ is bounded from above on $\D^n$;
\item
$\Psi$ is circled, i.e. $\Psi (z_1,\dots,z_n)$ depends only on $(|z_1|,
\dots, |z_n|)$;
\item
for any $c>0$, $\Psi (|z_1|^c, \dots, |z_n|^c) =c \Psi (|z_1|,
\dots, |z_n|)$.
\end{enumerate}

As a consequence,
$(dd^c \Psi)^n = \tau_\Psi \delta_0$ for some $\tau_\Psi \ge 0$.
\end{defn}

Notice that if $\Psi_1 \in PSH(\D^n)$, $\Psi_2 \in PSH(\D^m)$, 
and they are both local indicators, then
$$
\Psi (z,z'):= \max ( \Psi_1(z), \Psi_2(z'))
$$
defines a local indicator on $\D^{n+m}$. 

We need to remove the restriction to a single coordinate system in Definition \ref{LRindic}.

\begin{defn}
\label{genindic}
We call $\Psi$ a \emph{generalized local
indicator},
and we write $\Psi \in \mathcal I$ if there exists $U$
a neighborhood of $0$,
$\Psi_0 \in \mathcal I_0$
and a  one-to-one linear map $L$ of
$\C^n$ to itself such that $L(U ) \subset \D^n$
and $\Psi = \Psi_0 \circ L$.
\end{defn}

We will concentrate on a class of simple examples. Given two vectors $z,
w \in \C^n$, their standard hermitian product is denoted by $z\cdot \bar
w := \sum_j z_j \bar w_j$. We also write $\| z \| := |z \cdot \bar z |^{1/2}.$

\begin{defn}
\label{elemindic}
We say that $\Psi$ is an \emph{elementary local indicator} if there
exists a basis
$\{v_1, \dots, v_n\}$ of vectors of $\C^n$ and scalars $m_j \in \R_+$,
$1\le j \le n$,
such that for $z\in \D^n$,
\begin{equation}
\label{elindform}
\Psi(z) = \max_ {1\le j \le n} m_j \log |z \cdot \bar v_j | .
\end{equation}
\end{defn}

One easily checks that any elementary local indicator is a generalized local indicator.
The most interesting case is the one for which the basis
is orthornormal. In fact, it is essentially the only case.

\begin{lemma}
\label{ortho}
Given  an elementary local indicator $\Psi$ as in Definition \ref{elemindic} there exists an
orthonormal basis
$\{\tilde v_1, \dots, \tilde v_n\}$ of $\C^n$ such that the associated elementary local
indicator $\tilde \Psi (z) :=
\max_ {1\le j \le n} m_j \log |z \cdot \overline {{\tilde v}_j} |$
verifies $\tilde \Psi - \Psi \in L^\infty(\D^n)$.
\end{lemma}

As a consequence, we could have restricted the map $L$ in Definition
\ref{genindic} to
be unitary, and it would not have changed things in any essential way.

The proof of Lemma \ref{ortho} is given in Section \ref{easypf} below.

\begin{lemma}
 \cite[example in Section
3]{lelongrash}, \cite{Rash}
\label{mass}
If $\Psi$ is an elementary local indicator,
then $\tau_\Psi = m_1 \cdots m_n.$
\end{lemma}

We take the same definition of the generalized Green function as in
\cite{lelongrash}.

\begin{defn}
\label{defgreen}
Let $\Omega $ be a bounded domain in $\C^n$. 
Given 
$$
S := \{ (a_j,\Psi_j), 1\le j \le N\}, \mbox{ where } a_j \in \Omega, 
a_j \neq a_k \mbox{ for }j \neq k, 
\Psi_j \in \mathcal I,
$$
its \emph{Green function} is 
\begin{multline*}
G_S  :=
\sup  \left\lbrace u \in PSH(\Omega, \mathbb R_-) :  u(z) \le
\Psi_j (z) + C_j,
\right.
\\
\left.
\mbox{ for } z \mbox{ in a neighborhood of } a_j, j=1,...,N \right\rbrace .
\end{multline*}
\end{defn}

To generalize the Lempert function, the first step is to quantify the way in which an
analytic disk, i.e. an element of $Hol (\D,\Omega)$, meets a
pole provided with a generalized local indicator.

\begin{defn}
\label{multphia}
Let $\alpha \in \D$, $a \in \Omega$, $\Psi \in \mathcal I$. Then the
\emph{multiplicity} of $\varphi \in  Hol (\D,\Omega)$
at $\alpha$, with respect to $a$, is given by
\begin{eqnarray*}
\mbox{If } \varphi (\alpha) =a, & \mbox{ then } &
m_{\varphi, a, \Psi} (\alpha) :=
\min \left( \tau_\Psi , \liminf_{\zeta \to 0}
\frac{\Psi(\varphi(\alpha+\zeta)-a)}{\log |\zeta|} \right);
\\
\mbox{if } \varphi(\alpha) \neq a , &\mbox{ then }&
m_{\varphi, a, \Psi} (\alpha) :=0.
\end{eqnarray*}
\end{defn}

Notice that if $\Psi_1-\Psi_2$ is locally bounded near the origin, then
$m_{\varphi, a, \Psi_1} (\alpha)
=m_{\varphi, a, \Psi_2} (\alpha)$.

The quantity $\liminf_{\zeta \to 0}
\frac{\Psi(\varphi(\alpha+\zeta)-a)}{\log |\zeta|}$ is exactly the Lelong number at $0$
of the subharmonic function $\Psi \circ \varphi$, compare with \cite[pp. 334--335]{RashSig}.
Truncating at the level of the local Monge-Amp\`ere mass $\tau_\Psi$ will turn out
to be convenient in Definition \ref{defLempert}, and the proofs that use it. 

It is useful to see what this means in the case of elementary local indicators.

{\bf Elementary examples.}

Suppose that $\alpha=0$, $a=0$, and that $\Psi(z) = \max_ {1\le j \le n}
m_j \log |z_j | .$ We write
$$
\varphi(\zeta) = (\varphi_1(\zeta), \dots , \varphi_n(\zeta)),
$$
and define the valuations
$$
\nu_j := \nu_j (0,\varphi) := \min \{ k : (\frac{d}{d\zeta})^k
\varphi_j(0) \neq 0 \}.
$$
Then we have
\begin{equation}
\label{multex}
m_{\varphi, 0, \Psi} (0) = \min \left( \min_ {1\le j \le n} m_j \nu_j ,
\prod_{j=1}^n m_j \right).
\end{equation}

{\bf Example 1.}

If $m_j=1$ for all $j$, $m_{\varphi, 0, \Psi} (0)=1$ if
$\varphi(0)=0$, $m_{\varphi, 0, \Psi} (0)=0$
otherwise. This is the basic case where one just records whether a point
has been hit by the analytic disk or not.

{\bf Example 2.}

In more general cases, the use of an elementary local indicator will impose higher-order
differential conditions on the
map $\varphi$. For instance, if  $m_1=2$ and $m_j=1$, $2\le j \le n$, then
\begin{eqnarray*}
m_{\varphi, 0, \Psi} (0)& =& 0 \mbox{ if } \varphi(0) \neq 0;
\\
m_{\varphi, 0, \Psi} (0)& =& 1 \mbox{ if } \varphi(0) = 0 \mbox{ and
} \varphi_j'(0) \neq 0 \mbox{ for some }j \in \{2,\dots, n\};
\\
m_{\varphi, 0, \Psi} (0)& =& 2 \mbox{ if } \varphi(0) = 0 \mbox{ and
} \varphi_j'(0) = 0\mbox{ for any }j \in \{2,\dots, n\}.
\end{eqnarray*}

\begin{defn}
\label{defLempert}
Given a system $S$ as
in Definition \ref{defgreen},
we write $\tau_j := \tau_{\Psi_j}$. 

Let $\varphi \in
Hol(\D,\Omega)$
and $A_j \subset \D$, $1\le j \le N$.
We say that $(\varphi, (A_j)_{1\le j \le N})$ is \emph{admissible (for $S$, $z$)} if 
$$
\varphi (0)=z;
\quad
A_j \subset \varphi^{-1}(a_j)
\mbox{ and }
\sum_{\alpha \in A_j}
m_{\varphi, a_j, \Psi_j} (\alpha) \le \tau_j , 1\le j \le N.
$$
In this case, we write (with the convention that $0 \cdot \infty =0$)
$$
\mathcal S (\varphi, (A_j)_{1\le j \le N}):= \sum_{j=1}^N \sum_{\alpha
\in A_j} m_{\varphi, a_j, \Psi_j} (\alpha)
\log |\alpha|  .
$$
Then
the generalized Lempert function is defined by
\begin{multline*}
\mathcal L^\Omega_S (z):=
\mathcal L_S (z)
\\
:= \inf
\left\lbrace
\mathcal S (\varphi, (A_j)_{1\le j \le N}) : (\varphi, (A_j)_{1\le j \le N})
\mbox{ is admissible for }S, z 
\right\rbrace .
\end{multline*}
\end{defn}

Notice that we allow any of the $A_j$ to be the empty set (in which case
the $j$-th term drops from the sum). 

Consider the \emph{single poles case} where 
\begin{equation}
\label{sglp}
\mbox{ for each }j, \quad \Psi_j (z) =  \max_ {1\le l
\le n}  \log |z_l | ,\mbox{ or }\Psi_j (z) =  \log \|z \|
\end{equation}
-- it is the same, since
both functions differ by a bounded term near $0$; in fact, one could 
use any norm that is homogenous under complex scalar multiplication. 

In this case, $\tau_j=1$ for every $j$. With a slight abuse of notation, we write $S = \{ a_1, \dots , a_N\}.$
Then  $\mathcal L_S (z)= \min_{S'\subset S}
\ell_{S'}(z)$, where $\ell_S$ is defined in \eqref{coman}.
And in fact $\min_{S'\subset S}
\ell_{S'}(z) = \ell_{S}(z)$ \cite{NiPfl} (see also \cite{Wikstrom}, \cite{WikstromAMS}
for the case when the domain $\Omega$ is convex).

The Lempert function is different from the functionals considered by
Poletsky and others
in that it is restricted to one pre-image per pole $a_j$ (thus the Lempert function
can fail to be equal to the corresponding Green function). In our
definition, the number of pre-images
per pole is bounded above by the Monge-Amp\`ere mass at that pole of
its generalized local indicator. In \cite{TraoTh}, each pole only could have one
pre-image, but (essentially) $\varphi$ had
to hit the pole with maximum multiplicity at that pre-image.

Although Definition \ref{defLempert} may seem contrived, it is required to obtain
the reasonable convergence theorem \ref{conv}. See the discussion in Section \ref{compprev}.

We remark right away that the usual relationship holds between this
generalized Lempert function and
the corresponding Green function.

\begin{lemma}
\label{inegGL}
For $\Omega $ a bounded domain, 
for any system $S$ as in Definition \ref{defgreen}, for any $z \in \Omega$, $G_S (z) \le \mathcal L_S (z)$.
\end{lemma}

\begin{proof}
If $\varphi \in Hol(\D,\Omega)$, and $u \in PSH_- (\Omega)$ is a member
of the defining family for the Green
function of $S$, then $u \circ \varphi$ is subharmonic and negative on
$\D$. Furthermore, if $(\varphi, (A_j)_{1\le j \le N})$ is admissible (for $S$, $z$)
and $\alpha \in A_j$, then given any $\eps>0$, for $|\zeta|$ small enough,
$$
u \circ \varphi (\alpha +\zeta ) \le
C_j + \Psi_j (\varphi (\alpha +\zeta )-a_j)
\le
C_j +  (m_{\varphi, a_j, \Psi_j} (\alpha) - \eps ) \log |\zeta| .
$$
So $u \circ \varphi$ is a member of the defining family for the Green
function
on $\D$ with poles $\alpha$ and weights $m_{\varphi, a_j, \Psi_j}
(\alpha) - \eps$
at $\alpha$. This implies that
$$
u \circ \varphi (\zeta) \le
\sum_{j=1}^N \sum_{\alpha \in A_j} (m_{\varphi, a_j, \Psi_j} (\alpha) -
\eps )
\log \left| \frac{\alpha -\zeta}{1-\zeta \bar \alpha} \right| .
$$
Letting $\eps$ tend to $0$ and setting $\zeta=0$,
we get $u (z) \le \mathcal S (\varphi, (A_j)_{1\le j \le N})$.

Passing to the supremum over $u$, 
then to the infimum
over $(\varphi, (A_j)_{1\le j \le N})$, we get the Lemma.
\end{proof}

\section{Main Results}
\label{main}

We start with a remark.

\begin{lemma}
\label{smallerS}
If $S$ is as in Definition \ref{defgreen}, $1\le N'\le N$, and
$$
S':=  \{ (a_j,\Psi_j), 1\le j \le N'\},
$$
then for any $z \in \Omega$, $\mathcal L_{S'} (z) \ge \mathcal L_{S} (z)$.
\end{lemma}

\begin{proof} 
If we take $A_j=\emptyset$ for $N'+1 \le j \le N$, any member of
the defining family for $\mathcal L_{S'} (z)$
becomes a member of the defining family for $\mathcal L_{S} (z)$, and
the sum remains the same.
\end{proof}

The above lemma goes in the direction of monotonicity of the Lempert
function with respect to its system of poles.
For the Green function, it is immediate that
the more poles there are, the more negative the function must be. 
More generally the more negative the
generalized local indicators are (removing a pole 
corresponds to replacing a local  indicator by $0$), the more negative the function must be. 
This is not immediately
apparent in 
Definition \ref{defLempert}, but it does hold 
for elementary local indicators.

\begin{theorem}
\label{monotone}
Let $\Omega$ be  a bounded domain in $\C^n$,
$$
S := \{ (a_j,\Psi_j), 1\le j \le N\},
S' := \{ (a_j,\Psi'_j), 1\le j \le N\},
\mbox{ where } a_j \in \Omega,
$$
and $\Psi_j$, $\Psi'_j$, are elementary local indicators such that $\Psi_j \le
\Psi'_j + C_j$ in a neighborhood of $0$, $C_j \in \R$, $1\le j \le N$.
Then $\mathcal L_{S'} (z) \ge \mathcal L_{S} (z)$, for all $z \in \Omega$.
\end{theorem}

The proof is given in Section \ref{pfmono}.

Now we turn to a result about the convergence of some families of
(ordinary) Lempert functions with single poles, whose limits can be
described naturally as generalized Lempert functions.
Note that the proof of
this next theorem doesn't require the relatively difficult Theorem
\ref{monotone}, only the easy Lemma \ref{smallerS}.

For $z\in \C^n \setminus \{0\}$, we denote by $[z]$ the equivalence
class of $z$ in the complex projective space $\mathbb P^{n-1}$.

\begin{theorem}
\label{conv}
Let $\Omega$ be  a bounded and convex domain in $\C^n$.
Let $0\le M \le N$ be integers.
For $\eps$ belonging to a neighborhood of $0$ in $\C$, 
using the simplified notation of the single pole case \eqref{sglp}, let
$$
S(\eps) := \left\lbrace
a_j(\eps), 1\le j\le M ; a'_j(\eps), a''_j(\eps), M+1\le j\le N
\right\rbrace
\subset \Omega.
$$
Suppose that all the points of $S(\eps)$ are distinct for any fixed $\eps$,
that 
$$
\lim_{\eps\to 0} a_j(\eps) = a_j \in \Omega, 1\le j\le M ;
$$
$$
\lim_{\eps\to 0} a'_j(\eps) = \lim_{\eps\to 0} a''_j(\eps) =a_j \in
\Omega, M+1\le j\le N ;
$$
and that
\begin{equation}
\label{limproj}
\lim_{\eps\to 0} [a''_j(\eps)-a'_j(\eps)] = [v_j],
\end{equation}
where the limit is with respect to the
distance in $\mathbb P^{n-1}$ and the representative $v_j$ is chosen of unit norm.
Let $\Psi_j (z) := \log \| z \|$, $1\le j\le M$. Denote by $\pi_j$ the
orthogonal projection
onto $\{v_j\}^\bot$, $M+1\le j\le N $, and by $\Psi_j$ the generalized
local indicator
$$
\Psi_j (z) := \max (\log \| \pi_j(z)\|, 2 \log |z\cdot \bar v_j|), \quad
M+1\le j\le N.
$$
Set $S:=\{(a_j,\Psi_j), 1 \le j \le N\}$. Then
$$
\lim_{\eps\to 0} \ell_{S(\eps)} (z) =
\lim_{\eps\to 0} \mathcal L_{S(\eps)} (z) = \mathcal L_{S} (z) \mbox{
for all } z \in \Omega.
$$
\end{theorem}

Remarks : (a) as in the comments after \eqref{sglp}, one could replace
$\Psi_j$ by an elementary local indicator; (b) the convexity requirement
is imposed by Lemma \ref{relax}, and we conjecture that it is not essential.

Note that in the case where $a'_j(\eps)=a_j$ does not depend on $\eps$,
the hypothesis \eqref{limproj}
means that the point $a''_j(\eps)$ converges to a limit
in the blow-up of $\C^n$ around the point $a_j$.

It seems to us that this is the only reasonable convergence
result that can be obtained for a family of ordinary Lempert functions. If  \eqref{limproj}
is not satisfied, one can find two distinct
limit points for our family of Lempert functions. Thus hypothesis  \eqref{limproj}
is required.

We are restricting ourselves
to the case where no more than two points converge to the same
point: examples where three points converge to 
the origin in the bidisk are explicitly studied in \cite{Th3p}, and show that the
situation leads to results that probably can't be described in terms of our
generalized local indicators. 

The proof is given in Section \ref{pfconv}.

\section{Proof of Lemma \ref{ortho}}
\label{easypf}

Multiplying one of the vectors $v_j$ by a scalar only
modifies the function $\Psi$ by a bounded additive term, 
so it will be enough to
exhibit an orthogonal basis of
vectors complying with the conclusion of the Lemma.

Renumber the vectors $v_j$ so that we have $0\le m_1 \le \cdots \le
m_n$. Using the Gram-Schmidt orthogonalization
process, we produce an orthogonal system of vectors $\tilde v_k$ such
that 
$\mbox{Span} (\tilde v_1, \dots, \tilde v_k) = \mbox{Span} ( v_1,
\dots,  v_k)$ for any $k$, $1\le k \le n$.

We proceed by induction on the dimension $n$. When $n=1$ the
property is immediate. Assume that the result
holds up to dimension $n-1$. Write
$$
\Psi_1(z) := \max_ {1\le j \le n-1} m_j \log |z \cdot \bar v_j | , \quad
\tilde \Psi_1 (z) :=
\max_ {1\le j \le n-1} m_j \log |z \cdot \overline {{\tilde v}_j}  |.
$$
Denote $z_n:= z \cdot \overline {{\tilde v}_n} $.

It is enough to obtain the estimates on a neighborhood $U$ of $0$.
We choose it so that for $z \in U$, $|z_n|\le 0$, $ \Psi_1 (z), \tilde \Psi_1 (z) \le 0$. 
Since $v_n=\tilde v_n - w$, where $w \in \mbox{Span}( v_1, \dots,
v_{n-1})$,
we have
\begin{multline}
\label{psiw}
\Psi (z) =
\max (  \Psi_1 (z'), m_n \log |z_n-z'\cdot \bar w |), \\
\tilde \Psi (z) =
\max ( \tilde \Psi_1 (z'), m_n \log |z_n |),
\end{multline}
where $z'$ stands for the orthogonal projection of $z$ on $\mbox{Span}(
v_1, \dots,  v_{n-1})$
$=\mbox{Span}( \tilde v_1, \dots,  \tilde v_{n-1})$. By the induction
hypothesis, $ \Psi_1=\tilde \Psi_1 + O(1),$ so it is enough to prove that 
$$
\Psi'(z):= \max (  \Psi_1 (z'), m_n \log |z_n |)
$$
differs from $\Psi(z)$ by a bounded additive term.

There is a constant $C_0>0$ such that $\Psi_1 (z') \ge m_{n-1}
\log \| z' \| - \log C_0$, for 
$z'\in U$. Choose a constant $A>1$ large
enough so that $ \|w\| (C_0/A)^{1/ m_{n-1}} < 1/2$.

Then, since $ \Psi_1 (z) \le 0$ and $m_{n-1}\le m_{n}$,
\begin{multline}
\label{zw}
|z'\cdot \bar w |  \le \|w\| C_0^{1/ m_{n-1}} \exp(
\frac{\Psi_1 (z') }{m_{n-1}})
\\
\le \|w\| C_0^{1/ m_{n-1}} \exp(
\frac{\Psi_1 (z') }{m_{n}}).
\end{multline}

{\it Case 1.} $\Psi_1 (z') \ge m_n \log |z_n| - \log A.$

By the inequality above, $\Psi'(z) \le \Psi_1 (z') + \log A \le  \Psi (z) + \log A$.

On the other hand,
using \eqref{zw}, we get
$$
 |z_n-z'\cdot \bar w |^{m_{n}} \le  
\left(A^{1/m_{n-1}} +  \|w\| C_0^{1/m_{n-1}} \right)^{m_{n}}  \exp( \Psi_1 (z') ),
$$
so $\Psi(z) \le \Psi_1 (z') +O(1) \le \Psi' (z) +O(1)$.

{\it Case 2.} $\Psi_1 (z') \le m_n \log |z_n| - \log A.$

Then \eqref{zw} and the choice of $A$ imply
$$
|z'\cdot \bar w | \le \|w\| C_0^{1/ m_{n-1}} \exp\left( \log |z_n| -\frac{ \log A}{m_{n-1}} \right)
\le \frac12 |z_n|,
$$
thus \eqref{psiw} implies that
$$
\Psi'(z) +\log \frac12 \le \Psi(z) \le \Psi'(z) +\log \frac32.
$$

\section{Proof of Theorem \ref{monotone}}
\label{pfmono}


Without loss of generality, we may assume that $\tau'_j>0$ for all $j$.

We have
$\Psi_j \le \Psi'_j + C_j$ in a neighborhood of $0$ and
$$
\mbox{supp}\, (dd^c \Psi_j)^n \subset \{ 0 \}, \quad 
\mbox{supp}\, (dd^c \Psi'_j)^n \subset \{ 0 \}.
$$
Thus it follows from
Bedford and Taylor's comparison theorem \cite{BT}, \cite[p. 126, 
Theorem 3.7.1]{Kli} that 
$\tau_j \ge \tau'_j>0$. For any $\alpha, a_j$, 
\begin{equation}
\label{mgtm}
m_{\varphi, a_j,
\Psi_j} (\alpha) \ge
m_{\varphi, a_j, \Psi'_j} (\alpha).
\end{equation}
Therefore 
\begin{equation}
\label{comp_sums}
\sum_{j=1}^N \sum_{\alpha \in A_j} m_{\varphi, a_j, \Psi_j} (\alpha)
\log |\alpha|
\le
\sum_{j=1}^N \sum_{\alpha \in A_j} m_{\varphi, a_j, \Psi'_j} (\alpha)
\log |\alpha| .
\end{equation}
To finish the proof, it suffices to show that the family
over which we take the infimum
is smaller for $\mathcal L_{S'} (z)$ than the one for $\mathcal L_{S} (z)$. 
This can be checked
for each $j$ separately, hence we drop the index $j$.

\begin{lemma}
\label{comp_sets}
Let $\Omega$ be a bounded domain in $\C^n$. 
If $\Psi, \Psi'$ are elementary local indicators such that $\Psi \le
\Psi' + C$ and
$\tau':=\tau_{\Psi'}>0$,
if $A\subset \D$, $a\in \Omega$ and $\varphi \in Hol(\D,\Omega)$ verify
$$
 \sum_{\alpha \in A} m_{\varphi, a, \Psi'} (\alpha) \le \tau',
$$
then
$$
 \sum_{\alpha \in A} m_{\varphi, a, \Psi} (\alpha) \le \tau:=\tau_{\Psi}.
$$
\end{lemma}

\begin{proof}
Since the point $a$ plays no role, we assume
$a=0$ and write $m_{\varphi, 0, \Psi} (\alpha)=
m_{\varphi,  \Psi} (\alpha)$. By \eqref{mgtm}, we may assume that 
 $ m_{\varphi,  \Psi} (\alpha) >0$ and the sums in \eqref{comp_sums} will not change. 

Using Lemma \ref{ortho}, we reduce ourselves to the case where the
elementary local indicators are given by orthonormal systems of
vectors. We use the same ``valuations" as in the Elementary Example:
$$
\nu_j (\alpha) := \nu_j (\alpha ,\varphi) := \min \{ k :
(\frac{d}{d\zeta})^k (\varphi (\zeta) \cdot \bar v_j ) (\alpha )\neq 0 \},
$$
and $\nu'_j (\alpha)$ is defined analogously using the vectors $v'_j$.

{\bf Case 1.}
There exists $\alpha_0$ such that $m_{\varphi,  \Psi'} (\alpha_0) = \tau'$.

Then the hypothesis of Lemma \ref{comp_sets} implies that for all
$\alpha \in A\setminus \{\alpha_0\}$, $m_{\varphi,  \Psi'} (\alpha) =0$,
so that $\min_{1\le k \le n} m'_k \nu'_k(\alpha)=0$.
Since $\tau'>0$, we have $m'_k>0$ for all $k$, so there must exist
 $k$ such that $\nu'_k(\alpha)=0$. Then $\varphi(\alpha)\neq 0$, which
implies
that $m_{\varphi,  \Psi} (\alpha) =0$, and
$$
\sum_{\alpha \in A} m_{\varphi,  \Psi} (\alpha) =
 m_{\varphi,  \Psi} (\alpha_0) \le \tau,
$$
by definition of the multiplicity.

{\bf Case 2.}
For all $\alpha \in A$, $m_{\varphi,  \Psi'} (\alpha) < \tau'$.

Therefore $m_{\varphi,  \Psi'} (\alpha) = \min_{1\le k \le n} m'_k
\nu'_k(\alpha)$, and
since we always have $m_{\varphi,  \Psi} (\alpha) \le \min_{1\le k \le
n} m_k \nu_k(\alpha)$, it becomes enough to work with those quantities
in \eqref{comp_sums}. By dividing by $\tau$ and $\tau'$
respectively, it will be enough to prove the following Lemma.
\end{proof}

\begin{lemma}
\label{dirmult}
Under the hypotheses of Lemma \ref{comp_sets} and Case 2 above, for each $\alpha \in A$,
$$
\frac{\min_{1\le k \le n} m'_k \nu'_k(\alpha)}{\prod_{k=1}^n m'_k}
\ge
\frac{\min_{1\le k \le n} m_k \nu_k(\alpha)}{\prod_{k=1}^n m_k} .
$$
\end{lemma}

\begin{proof}
Since we are now dealing with a single $\alpha$, we also drop it from
the notation.

We introduce a binary relation on the
index set $\{1,\dots,n\}$.

\begin{defn}
Given $k, l \in \{1,\dots,n\}$, we say that $k \mathcal R l$ if and only if
$v_k \cdot \bar v'_l \neq 0$.
\end{defn}

\begin{lemma}
\label{relmult}
If $\Psi' + C \ge \Psi$ and $k \mathcal R l$, then $m_k \ge m'_l$.
\end{lemma}

\begin{proof}
For any nonzero $\lambda \in \mathbb C$,
$$
\Psi'(\lambda v'_l)= m'_l \log |\lambda| + m'_l \log \|v'_l\|^2,
$$
while, for $|\lambda|$ small enough,
$$
\Psi(\lambda v'_l)= \max_{1\le j \le n}
\left( m_j (\log |\lambda| +  \log |v_j \cdot \bar v'_l |) \right)
=(\min_{k: k \mathcal R l} m_k) \log |\lambda| + O(1),
$$
therefore by letting $\lambda$ tend to $0$ we see that
$\min_{k: k \mathcal R l} m_k \ge m'_l$.
\end{proof}

\begin{lemma}
\label{relval}
If $\Psi' + C \ge \Psi$, then
\begin{enumerate}
\item
$\nu'_l \ge \min \{ \nu_k : k \mathcal R l \},$
\item
$\nu_k \ge \min \{ \nu'_l : k \mathcal R l \}.$
\end{enumerate}
\end{lemma}

\begin{proof}
We will use and prove part (1) only. The
other one has a similar proof.

Since $v'_l$ is orthogonal to $v_k$ unless $k \mathcal R l$, we must
have complex scalars $c_k$ such that $v'_l =\sum_{k : k \mathcal R l}
c_k v_k$, thus for $\varphi $ as in Lemma \ref{comp_sets},
$$
\varphi(\zeta) \cdot \bar v'_l =
\sum_{k : k \mathcal R l} \bar c_k \varphi(\zeta) \cdot \bar v_k .
$$
Now take $m < \nu_k = \nu_k(\alpha, \varphi)$, for all $k$ such that $k \mathcal R l$. Then
$$
(\frac{d}{d\zeta})^m (\varphi \cdot \bar v'_l ) (\alpha )
= \sum_{k : k \mathcal R l} \bar c_k (\frac{d}{d\zeta})^m (\varphi
(\zeta) \cdot \bar v_k ) (\alpha) = 0,
$$
so we must have $\nu'_l > m$, which proves the result.
\end{proof}

Now renumber the vectors $v'_l$ so that $\min_k (m'_k \nu'_k) = m'_1
\nu'_1$.
Pick an index $k_0$ such that $k_0 \mathcal R 1$ and $\nu_{k_0} =
\min \{ \nu_k : k \mathcal R 1 \}$. By renumbering the vectors $v_k$,
we may assume $k_0=1$. By Lemma \ref{relval}, we may assume 
$\nu'_1 \ge \nu_1$.

The conclusion of Lemma \ref{dirmult} thus reduces to:
\begin{equation}
\label{ineqfin}
\frac{ \nu'_1}{\prod_{k=2}^n m'_k}
\ge
\frac{\nu_1}{\prod_{k=2}^n m_k} .
\end{equation}

This is a consequence of the next result. 

\begin{lemma}
\label{order}
There exists a bijection $\sigma$ from $\{2,\dots,n\}$ onto itself
such that for any $l \in \{2,\dots,n\}$, $ \sigma(l)\mathcal R l $.
\end{lemma}

This Lemma will be proved below. It implies
$$
\prod_{k=2}^n m_k = \prod_{l=2}^n m_{\sigma(l)} \ge \prod_{l=2}^n m'_l,
$$
by Lemma \ref{relmult}, so \eqref{ineqfin} holds and this concludes
the proof of Lemma \ref{dirmult}.
\end{proof}

\begin{proof*}{\it Proof of Lemma \ref{order} }

Denote $A:=(a_{kl})_{2\le k,l \le n} := (v_k \cdot \overline{v'_l})_{2\le k,l \le n}$.
First we prove that this matrix is non singular. Let $\pi$ be the orthogonal projection
on $\{v'_1\}^\bot$. If $\mbox{rank}\, \{\pi(v_k), 2 \le k \le n\} < n-1$, there exists
$w \in \{v'_1\}^\bot$, $w\neq 0$, such that $w \bot
\pi(v_k)$, $2\le k \le n$. This implies $w \bot v_k$, $2\le k \le n$. 
Since we have orthogonal bases, $v_1= \lambda w$,
for some $\lambda \in \C$. So $v_1 \cdot
\overline{v'_1}=0$, which contradicts the fact that $1\mathcal R 1$. 

We construct the bijection $\sigma$ by induction  on $n$. For $n=2$ it's obvious. 
Suppose that the
property holds for $n-1$. Then
$$
0\neq \mbox{\rm det} A = \sum_{k=2}^n (-1)^k a_{k2} \mbox{\rm det} A_k,
$$
where $A_k$ stands for the minor matrix with the first column and the
$k$-th row
removed. There must be some $k$ for which $a_{k2} \mbox{\rm det} A_k
\neq 0$. Let
$\sigma(2)=k$; the induction hypothesis gives us a bijection $\sigma'$ from
$\{3,\dots,n\}$ to $\{2,\dots,n\}\setminus \{k\}$ such that
$a_{\sigma'(l)l}\neq 0$,
and this finishes the proof.
\end{proof*}

\section{Proof of Theorem \ref{conv}}
\label{pfconv}

\begin{proof*}

First observe that we can relax 
the conditions
used in Definition \ref{defLempert}.
\begin{lemma}
\label{relax}
Let $\Omega$ be a convex bounded domain in $\mathbb C^n$
containing the origin, and let $z\in \Omega$.

(i) Let $a_j \in \Omega$, $\Psi_j \in \mathcal I$ and, as in Definition \ref{defgreen}
$$
S := \left\lbrace
(a_j, \Psi_j), 1\le j\le N
\right\rbrace .
$$

Suppose that for any $\delta >0$, there exists a map $\varphi^\delta $
holomorphic from
$\D$ to $(1+\delta)\Omega$ and sets $(A_j(\delta))_{1\le j \le N}$ such
that
$(\varphi^\delta, (A_j(\delta))_{1\le j \le N})$ is admissible for $
S,z$ with respect to $(1+\delta)\Omega$ and
$$
\mathcal S (\varphi^\delta, (A_j(\delta))_{1\le j \le N}) \le \ell + h(\delta),
$$
where $ h(\delta)\ge 0$, $\lim_{\delta \to 0} h(\delta)= 0$.
Then $\mathcal L^\Omega_S (z) \le \ell.$

(ii) For $\eps$ in a neighborhood $V$ of $0$ in $\C$, let $a_j(\eps) \in \Omega$, $1\le j\le N$.
$$
S(\eps) := \left\lbrace
(a_j(\eps), \Psi_j), 1\le j\le N
\right\rbrace.
$$
Let $g: V \longrightarrow \R_+^*$ be 
such that $\lim_{\eps\to0}g(\eps)=0$.
Then
$$
\limsup_{\eps\to0} \mathcal L^\Omega_{S(\eps)} (z)
\le \limsup_{\eps\to0} \mathcal L^{(1+g(\eps))\Omega}_{S(\eps)} (z).
$$
\end{lemma}

\begin{proof}
Without loss of generality, we may assume $z=0$.

Let 
$$
\Omega_r := \left\{ \varphi (\zeta) : \varphi \in Hol(\D, \Omega), \varphi(0)=0, |\zeta|<r \right\}.
$$
A bounded convex domain is Kobayashi complete hyperbolic
\cite[Proposition 6.9 (b), p. 88]{Dineen}, so $\Omega_r $ is relatively compact in $\Omega $.
Let $\rho_\Omega$ stand for the Minkowski function of $\Omega$:
$$
\rho_\Omega(z):= \inf \{ r>0 : \frac{z}r  \in \Omega\}.
$$
We set $\gamma_\Omega (r) := \sup_{\Omega_r } \rho_\Omega.$
The function $\gamma$ is
increasing and continuous 
 from $(0,1)$ to itself.

For any $\mu \in (0,1)$, $\phi  \in Hol(\D, \C^n)$,
denote $\phi_\mu (\zeta) :=
\phi (\mu\zeta)$. Note that for any points and generalized local indicators,
$m_{\phi_\mu, a, \Psi} (\alpha/\mu) = m_{\phi, a, \Psi} (\alpha)$.

Take $\varphi^\delta$ as in Part (i) of the Lemma, in particular $\varphi^\delta(0)=0$, so by 
construction of $\gamma$, 
$$
\frac1{(1+\delta)}\varphi_\mu^\delta (\D)  \subset \gamma(\mu) \Omega .
$$
Choose some $\mu(\delta)$ such that $ \gamma(\mu(\delta)) =
(1+\delta)^{-1}$,
and set $\tilde  \varphi^\delta  := \varphi_{\mu(\delta)}^\delta $, then
$\tilde  \varphi^\delta \in Hol(\D, \Omega)$. Note that $\lim_{\delta\to
0} \mu(\delta)=1$, by the relative compactness of each $\Omega_r$.

Let
$$
\tilde A_j (\delta) := \left\lbrace \frac{\alpha}{\mu(\delta)}: \alpha
\in A_j (\delta), |\alpha|< \mu(\delta)
\right\rbrace .
$$
Then
\begin{multline}
\label{contract}
\left|
\mathcal S (\tilde \varphi^\delta, (\tilde A_j(\delta))_{1\le j \le N})
- \mathcal S (\varphi^\delta, (A_j(\delta))_{1\le j \le N})
\right|
\\
=
\left|
\sum_j \sum_{\alpha \in A_j, |\alpha|< \mu(\delta)} m_{\varphi^\delta,
a_j, \Psi_j} (\alpha) | \log \mu(\delta)|
-
\sum_j \sum_{\alpha \in A_j, |\alpha|\ge  \mu(\delta)}
m_{\varphi^\delta, a_j, \Psi_j} (\alpha) \log |\alpha|
\right|
\\
\le
2 (\sum_j \tau_{\Psi_j}) |\log \mu(\delta)|,
\end{multline}
and this last quantity tends to $0$, which concludes the proof of (i).

To prove (ii), take maps $\varphi^\eps$ and systems of points
$(A_j(\eps))$,
admissible for $S(\eps)$,
such that
$$
\lim_{\eps\to 0} \mathcal S ( \varphi^\eps, ( A_j(\eps))_{1\le j \le N}) =
\limsup_{\eps\to0} \mathcal L^{(1+g(\eps))\Omega}_{S(\eps)} (0).
$$
Use the above proof with $\delta=g(\eps)$ to construct
maps $\tilde \varphi^\eps$ into $\Omega$
and systems of points $(\tilde A_j(\eps))$, admissible for $S(\eps)$,
such that
\begin{equation*}
\left|
\mathcal S (\tilde \varphi^\eps, (\tilde A_j(\eps))_{1\le j \le N})
- \mathcal S (\varphi^\eps, (A_j(\eps))_{1\le j \le N})
\right|
\le
2 (\sum_j \tau_{\Psi_j}) |\log \mu(g(\eps))|,
\end{equation*}
and by definition $\mathcal S (\tilde \varphi^\eps, (\tilde
A_j(\eps))_{1\le j \le N}) \ge
\mathcal L^\Omega_{S(\eps)} (0) $.
\end{proof}

Consider as in Theorem \ref{conv} a bounded
convex domain $\Omega,$ and distinct points
$a_j \in \Omega$, $1\le j \le N$. Let $z \in \Omega \setminus \{a_j, 1 \le j \le N\}$ (otherwise
the property is trivially true).  Again we may assume $z=0$. 
By Lemma \ref{relax} applied to $S(\delta)=S$ for any $\delta$, 
 to show that
\begin{equation}
\label{liminfgr}
\mathcal L_S (z) \le \liminf_{\eps\to 0} \mathcal L_{S(\eps)} (z)=: \ell ,
\end{equation}
it will be
enough to provide: some increasing function $g$ such that $g(0)=0$ and,
 for any $\delta>0$, $\varphi^\delta \in Hol(\D,(1+g(\delta))\Omega)$
and subsets $(A_j^\delta )_{1\le j
\le N})$ of $\Omega$ such 
 $(\varphi^\delta , (A_j^\delta )_{1\le j
\le N})$ is
admissible for $S, z$ and that
$$
\mathcal S (\varphi^\delta, (A_j^\delta)_{1\le j \le N}) = \ell.
$$

The systems $S(\eps)$ all have single poles,
so the definition of $\ell$ means that there exist $\varphi_m \in Hol
(\D, \Omega)$,
$\eps_m\to 0$, and points
$\alpha_{j,m}, \alpha'_{j,m}, \alpha''_{j,m} \in \D $ such that
$\varphi_m ( \alpha_{j,m})= a_j(\eps_m),$ $1\le j \le M$, and
$\varphi_m ( \alpha'_{j,m})= a'_j(\eps_m),$ $\varphi_m (
\alpha''_{j,m})= a''_j(\eps_m),$ $M+1\le j \le N$; and they satisfy
$$
 \sum_{j=1}^M \log | \alpha_{j,m} | + \sum_{j=M+1}^N  
\left( \log |\alpha'_{j,m} | + \log | \alpha''_{j,m} | \right) = \ell
+ \delta (m) ,
$$
with $\lim_{m\to \infty} \delta (m) =0$.

Passing to a subsequence, for which we keep the same notations, we may
assume that
$\alpha_{j,m}\to \alpha_{j} \in \overline \D,$ $\alpha'_{j,m}\to
\alpha'_{j} \in \overline \D,$
$\alpha''_{j,m}\to \alpha''_{j} \in \overline \D$ as $m\to \infty$,
and that
$\varphi_m \to \tilde \varphi \in Hol (\D, \overline \Omega)$ uniformly on
compact subsets of
$\D$. Furthermore, by compactness of the unit circle, there exists
$\tilde {v_j} \in [v_j] \cap S^{2n-1}$ such that,
taking a further subsequence, 
$$
\lim_{m\to \infty}
\frac{a''_j(\eps_m)-a'_j(\eps_m)}{\| a''_j(\eps_m)-a'_j(\eps_m)
\|} = \tilde {v_j}.
$$

By renumbering the points and exchanging $a'_j$ and $a''_j$ as needed,
we may assume that there are integers $M' \le M$, $M\le N_1 \le  N_2 \le N_3
\le N$ such that
\begin{eqnarray*}
\alpha_j \in \D &\mbox{ for }& 1 \le j \le M' \\
\alpha_j \in \partial \D &\mbox{ for }& M'+1 \le j \le M\\
\alpha'_j = \alpha''_j\in \D &\mbox{ for }& M+1 \le j \le N_1 \\
|\alpha'_j| < |\alpha''_j| <1 &\mbox{ for }& N_1+1 \le j \le N_2 \\
|\alpha'_j |<1 , |\alpha''_j| =1 &\mbox{ for }& N_2+1 \le j \le N_3 \\
|\alpha'_j|=|\alpha''_j|=1 &\mbox{ for }& N_3+1 \le j \le N.
\end{eqnarray*}
Then
\begin{multline*}
\ell = \lim_{m\to\infty}
\left(
\sum_{j=1}^M \log |\alpha_{j,m}| + \sum_{j=M+1}^N \left(  \log |\alpha'_{j,m}| +
\log |\alpha''_{j,m}|\right) 
\right) \\
= \sum_{j=1}^{M'} \log |\alpha_{j}| + \sum_{j=M+1}^{N_1} 2 \log
|\alpha'_{j}|
+ \sum_{j=N_1+1}^{N_2} \left(  \log |\alpha'_{j}| +  \log |\alpha''_{j}| \right) 
+ \sum_{j=N_2+1}^{N_3}  \log |\alpha'_{j}| .
\end{multline*}

Now we choose
\begin{eqnarray*}
A_j = \{\alpha_j \} &\mbox{ for }& 1 \le j \le M' \\
A_j = \emptyset &\mbox{ for }& M'+1 \le j \le M\\
A_j = \{\alpha'_j \} &\mbox{ for }& M+1 \le j \le N_1 \\
A_j = \{\alpha'_j, \alpha''_j\} &\mbox{ for }& N_1+1 \le j \le N_2 \\
A_j = \{\alpha'_j \} &\mbox{ for }& N_2+1 \le j \le N_3 \\
A_j = \emptyset &\mbox{ for }& N_3+1 \le j \le N.
\end{eqnarray*}
Notice that $(\tilde \varphi, (A_j)_{1\le j \le N})$ hits the correct points
but doesn't necessarily
produce an admissible choice, because for some $j$, $N_1+1 \le j \le N_2$,
we could have
$$
m_{\tilde \varphi, a_j, \Psi_j} (\alpha'_j) + m_{\tilde \varphi, a_j,
\Psi_j} (\alpha''_j)
>  2=\tau_j.
$$

So, in order to apply Lemma \ref{relax} with $\delta\to 0$,
we set $A_j^\delta =A_j$ for
any $\delta>0$ and
$$
\tilde \varphi^\delta (\zeta) :=
\tilde \varphi (\zeta) + \delta \zeta
\left[ 
\prod_1^{M'} (\zeta-\alpha_j) \prod_{M+1}^{N_1} (\zeta-\alpha'_j)^2
\prod_{N_1+1}^{N_2} (\zeta-\alpha'_j) (\zeta-\alpha''_j)
\prod_{N_2+1}^{N_3} (\zeta-\alpha'_j) 
\right] \, v,
$$
where $v \in \C^n$ is a unit vector chosen such that $\pi_j(v) \neq 0$,
$N_1+1 \le j \le N_3$. For any $\alpha \in \cup_1^N A_j$,
$\tilde \varphi^\delta (\alpha) =
\tilde \varphi (\alpha)$. There is a constant $C>0$ such that
$\tilde \varphi^\delta  (\D) \subset \Omega + C \delta B(0,1)$.

All the following considerations apply when $\delta$ is small enough.

For $1\le j \le M'$, $m_{\tilde \varphi^\delta, a_j, \Psi_j}
(\alpha_j)=1$, because $\tilde \varphi^\delta$ takes on the correct
value, and the multiplicity cannot be more than $1=\tau_j$ in those cases.

For  $N_1+1\le j \le N_3$, we have
$$
\pi_j ((\tilde \varphi^\delta)'(\alpha'_j)) =
\pi_j ((\tilde \varphi)'(\alpha'_j)) + \delta p_j \pi_j(v),
$$
where $p_j$ is some complex scalar which doesn't depend on $\delta$, so
for $\delta >0$ and small enough, this projection doesn't vanish and we have
$m_{\tilde \varphi^\delta, a_j, \Psi_j} (\alpha'_j)=1$.  An analogous
reasoning shows that $m_{\tilde \varphi^\delta, a_j, \Psi_j} (\alpha''_j)=1$
for $N_1+1\le j \le N_2$.

For $M+1 \le j \le N_1$, we have
$$
(\tilde \varphi^\delta)'(\alpha'_j) = (\tilde \varphi)'(\alpha'_j),
$$
and by the uniform convergence on compact sets,
$$
(\tilde \varphi)'(\alpha'_j) =
\lim_{m\to\infty} \frac{\varphi_m (\alpha'_{j,m})-\varphi_m
(\alpha''_{j,m})}{\alpha'_{j,m}-\alpha''_{j,m}}
=
\lim_{m\to\infty}
\frac{a'_j(\eps_m)-a''_j(\eps_m)}{\alpha'_{j,m}-\alpha''_{j,m}},
$$
which must be colinear to $v_j$ by definition. Therefore
$m_{\tilde \varphi^\delta, a_j, \Psi_j} (\alpha'_j)=2$ for $M+1 \le j
\le N_1$.

Thus $(\tilde \varphi^\delta, (A_j)_{1\le j \le N})$
is admissible for $S, 0$ and $\mathcal S(\tilde \varphi^\delta, (A_j)_{1\le j \le N})=\ell$,
which proves \eqref{liminfgr}.

Now we need to show that
\begin{equation}
\label{limsupsm}
\mathcal L_S(z) \ge \limsup_{\eps\to0} \mathcal L_{S(\eps)}(z).
\end{equation}
We use Lemma \ref{relax}(ii). For any $\delta >0$,
we need to construct a positive function $g$ such that
$\lim_{\eps\to0}g(\eps)=0$
and, for $\eps$ small enough, $\varphi^\eps \in Hol(\D(1+g(\eps))\Omega)$
and sets $(A_j(\eps))_{1\le j \le N}$
such that
 $  (\varphi^\eps, (A_j(\eps))_{1\le j \le N})$ is
admissible for $S(\eps),0$
and
$$
\mathcal S (\varphi^\eps, (A_j(\eps))_{1\le j \le N}) \le \mathcal L_S(z) +\delta.
$$

We start
with an admissible choice $(\varphi, (A_j)_{1\le j \le N})$ for $S$,
such that
$$
\mathcal S (\varphi, (A_j)_{1\le j \le N}) \le \mathcal L_S (z) + \delta/2.
$$

To fix notations, suppose that, after renumbering and exchanging the points as needed,
there exist integers $M'\le M$, $N_1, N_2, N_3 \in \{ M, \dots N \}$
such that
\begin{eqnarray*}
A_j = \{\alpha_j \} &\mbox{ for }& 1 \le j \le M' ,\\
A_j = \emptyset &\mbox{ for }& M'+1 \le j \le M,\\
A_j = \{\alpha'_j \},
m_{\varphi, a_j, \Psi_j} (\alpha'_j)=2
 &\mbox{ for }& M+1 \le j \le N_1 ,\\
A_j = \{\alpha'_j, \alpha''_j\}, \alpha'_j \neq \alpha''_j &\mbox{ for }& N_1+1 \le j \le N_2 ,\\
A_j = \{\alpha'_j \},
m_{\varphi, a_j, \Psi_j} (\alpha'_j)=1 &\mbox{ for }& N_2+1 \le j \le N_3 ,\\
A_j = \emptyset &\mbox{ for }& N_3+1 \le j \le N.
\end{eqnarray*}

The definition of $\Psi_j$ (see the computations performed in the Elementary example)
implies that, for $M+1 \le j \le N_1$,
$\varphi'(\alpha'_j) \cdot \bar w = 0$, for any $w \in v_j^\perp$. We 
perturb $\varphi$ to make sure that, on the other hand,
$\varphi'(\alpha'_j) \cdot \bar v_j \neq 0$ in the same index range.
For $\eta(\eps) \in \mathbb C$ to be chosen later, set
\begin{multline*}
\tilde \varphi (\zeta) := \varphi (\zeta) +
\\
\eta(\eps) \,  \left[  \zeta \, \prod_1^{M'} (\zeta - \alpha_j)
\prod_{j=N_1+1}^{N_2} (\zeta - \alpha'_j)(\zeta - \alpha''_j)
\prod_{j=N_2+1}^{N_3} (\zeta - \alpha'_j) \right] 
\times
\\
\times
\left\{
\sum_{j=M+1}^{N_1} \left[ (\zeta - \alpha'_j) \prod_{M+1\le k \le N_1, k\neq j} (\zeta - \alpha'_k)^2 \right]
\, v_j
\right\}.
\end{multline*}
The map $\tilde \varphi$ depends on $\eps$ and is admissible again.

We have positive constants
$C_1, C_2, C_3$ such that 
\begin{itemize}
\item
$\tilde \varphi'(\alpha'_j)=\lambda_j v_j$, with $C_1^{-1} |\eta(\eps)| \le |\lambda_j | \le C_1 |\eta(\eps)| $,
\item
$\|\tilde \varphi - \varphi\|_\infty \le C_2 |\eta(\eps)|$, 
\item
$\tilde \varphi(\D) \subset (1+C_3|\eta(\eps)|) \Omega$;
\end{itemize}
in particular $\tilde \varphi$ will be bounded by constants
independent of $\eps$, along with all its derivatives on any given compact subset
of $\D$. 

For $M+1 \le j \le N$ and $\eps$ in a neighborhood of $0$, $a''_j(\eps) - a'_j(\eps) = n_j(\eps) v_j(\eps)$,
where $\|v_j(\eps)\|=1$, $\lim_{\eps\to 0} v_j(\eps) = v_j$ and $n_j(\eps) \in \C$.

For $|\eps|$ small enough, we now may define
\begin{multline*}
A_j(\eps) := A_j, \mbox{ for } 1 \le j \le M, N_1+1 \le j \le N, \\
\mbox{ and } A_j(\eps) :=\{ \alpha'_j, \alpha'_j +\frac{n_j(\eps)}{\lambda_j} \},
\mbox{ for } M+1 \le j \le N_1.
\end{multline*}
We shall need to add to 
$\tilde \varphi$ a vector-valued correcting term obtained by Lagrange
interpolation. To this end, we write
$B(\eps) := \cup_j A_j(\eps)$, and values to be interpolated,
$w (\alpha)$, for $\alpha \in B(\eps)$. Let
\begin{eqnarray*}
w(\alpha_j ) := a_j(\eps) - a_j = a_j(\eps) - \tilde \varphi(\alpha_j )
&\mbox{ for }& 1 \le j \le M' ,\\
w(\alpha'_j ) := a'_j(\eps) - a_j = a'_j(\eps) - \tilde \varphi(\alpha'_j )
&\mbox{ for }& M+1 \le j \le N_1 ,\\
w(\alpha'_j + \frac{n_j(\eps)}{\lambda_j}) := a''_j(\eps)
- \tilde \varphi(\alpha'_j + \frac{n_j(\eps)}{\lambda_j})
 &\mbox{ for }& M+1 \le j \le N_1 ,\\
w(\alpha'_j ) := a'_j(\eps) - a_j = a'_j(\eps) - \tilde \varphi(\alpha'_j )
 &\mbox{ for }& N_1+1 \le j \le N_2 ,\\
w(\alpha''_j ) := a''_j(\eps) - a_j = a''_j(\eps) - \tilde \varphi(\alpha''_j )
 &\mbox{ for }& N_1+1 \le j \le N_2 ,\\
w(\alpha'_j ) := a'_j(\eps) - a_j = a'_j(\eps) - \tilde \varphi(\alpha'_j )
 &\mbox{ for }& N_2+1 \le j \le N_3 .
\end{eqnarray*}
We denote by $P_\eps$ the solution to the interpolation problem
$$
\left(P(\alpha)= w(\alpha): \alpha \in B(\eps) \right).
$$
Let $\varphi^\eps := \tilde \varphi + P_\eps \in Hol(\D, \Omega^\eps)$. The
domain $\Omega^\eps$ will be specified below.
By construction
 $  (\varphi^\eps, (A_j(\eps))_{1\le j \le N})$ is
admissible for $S(\eps)$, and for $|\eps|$ small enough, 
$$
\mathcal S (\varphi^\eps, (A_j(\eps))_{1\le j \le N}) \le \mathcal L_S(z) +\delta,
$$
provided that, for $M+1 \le j \le N_1$,
\begin{equation}
\label{lambdacond}
\lim_{\eps\to 0}  \frac{n_j(\eps)}{\lambda_j} = 0,
\end{equation}

Now we need to show that the correction is small, more precisely
that we can choose $\eta(\eps)$ so that the above condition is satisfied
and $\lim_{\eps\to 0} \|P_\eps\|_\infty = 0$. Then we can choose a function
$g$ tending to $0$ such that 
$$
\Omega^\eps = (1+g(\eps) ) \Omega \supset (1+C_3|\eta(\eps)|) \Omega + B(0,\|P_\eps\|_\infty).
$$

Write $\Pi_\alpha$ for the unique (scalar) polynomial of degree
less or equal to $d:= \# B(\eps) - 1$
($d$ does not depend on $\eps$)
such that
$$
\Pi_\alpha (\alpha) = 1, \Pi_\alpha (\beta) = 0 \mbox{ for any }
\beta \in B(\eps) \setminus \{\alpha\}.
$$
Then
$$
P_\eps = \sum_{\alpha \in B(\eps)} \Pi_\alpha  w(\alpha).
$$
For $\alpha \in \bigcup_{1 \le j \le M, N_1+1 \le j \le N} A_j$,
 $\|\Pi_\alpha\|_\infty$ is uniformly bounded,
because $\mbox{dist} (\alpha, B(\eps) \setminus \{\alpha\}) \ge \gamma > 0$
with $\gamma$ independent of $\eps$. It also follows from the hypotheses of the
theorem and the choice of $w$ that
$$
\lim_{\eps\to 0} \, \max \{ \|w(\alpha)\| ,
\alpha \in \bigcup_{1 \le j \le M, N_1+1 \le j \le N} A_j\} = 0.
$$

For $M+1 \le j \le N_1$, we need an elementary lemma about Lagrange interpolation.
\begin{lemma}
\label{Lagrange}
Let $x_0, \dots , x_d \in \D$, $w_0, w_1 \in \C^n$. Suppose that there
exists $\gamma >0$ such that $|x_0-x_1|\le \gamma $ and
$\mbox{dist} ([x_0,x_1], \{x_2, \dots, x_d\})\ge 2 \gamma$,
where $[x_0,x_1]$ is the real line segment from $x_0$ to $x_1$.

Let $P$ be the unique
($\C^n$-valued) polynomial of degree less or equal to $d$ such that
$$
P(x_0)=w_0, P(x_1)=w_1, P(x_j)=0, 2 \le j \le d.
$$
Then there exist constants $L_1, L_0$ depending only on $\gamma$ and $d$ such
that
$$
\sup_{\zeta \in \D} \|P(\zeta)\| \le
L_1 \left\| \frac{w_1-w_0}{x_1-x_0} \right\| + L_0 \|w_0\|.
$$
\end{lemma}

We will prove this Lemma a little later.
It yields, for $M+1 \le j \le N_1$,
\begin{multline*}
\sup_{\zeta \in \D} \left\|
\Pi_{\alpha'_j}(\zeta) w(\alpha'_j )
+
\Pi_{\alpha'_j + \frac{n_j(\eps)}{\lambda_j}}  (\zeta)
w(\alpha'_j + \frac{n_j(\eps)}{\lambda_j})
\right\|
\\
\le
L_1  \left|\frac{\lambda_j} {n_j(\eps)}\right|
\left\|
 a''_j(\eps) - \tilde \varphi(\alpha'_j + \frac{n_j(\eps)}{\lambda_j})
\right\|
+
L_0 \|a'_j(\eps) - a_j\|.
\end{multline*}
We now estimate the first term in the last sum above.
By the Taylor formula,
$$
a''_j(\eps) - \tilde \varphi(\alpha'_j + \frac{n_j(\eps)}{\lambda_j})
=
a''_j(\eps) - a'_j(\eps) - n_j(\eps) v_j + R_2 (\eps)
=
n_j(\eps) ( v_j(\eps) -  v_j )  + R_2 (\eps),
$$
where $\| R_2 (\eps)\| \le C |n_j(\eps)|^2|\lambda_j|^{-2}$ with
$C$ a constant independent of $\eps$ by the boundedness of 
the derivatives of $\tilde \varphi$. Finally
\begin{multline*}
\sup_{\zeta \in \D} \left\|
\Pi_{\alpha'_j}(\zeta) w(\alpha'_j )
+
\Pi_{\alpha'_j + \frac{n_j(\eps)}{\lambda_j}}  (\zeta)
w(\alpha'_j + \frac{n_j(\eps)}{\lambda_j})
\right\|
\\
\le
C \left(
\|v_j(\eps) -  v_j\| |\eta(\eps)|
+
|n_j(\eps)||\eta(\eps)|^{-1}
+
\|a'_j(\eps) - a_j\|
\right).
\end{multline*}
To to satisfy \eqref{lambdacond}, we need to have $\lim_{\eps\to 0} n_j(\eps)/\eta(\eps)=0$;
to make sure, in addition, that the whole sum above tends to $0$ as $\eps$ tends to $0$,
 it will be enough to choose
$\eta(\eps)$ going to zero, but more slowly that $|n_j(\eps)|
=\|a''_j(\eps) - a'_j\|$, for $M+1 \le j \le N_1$.  
\end{proof*}

\begin{proof*}{\it Proof of Lemma \ref{Lagrange}}
Let
$$
Q(X,Y):= \prod_2^d \frac{X-x_k}{Y-x_k}.
$$
Then $Q$ and all of its derivatives are bounded for
$X \in \overline \D$ and $Y\in [x_0,x_1]$.
\begin{multline*}
P(X) = \frac{X-x_0}{x_1-x_0} Q(X,x_1) w_1 + \frac{X-x_1}{x_0-x_1} Q(X,x_0) w_0 \\
= \frac{w_1 - w_0}{x_1-x_0} (X-x_0) Q(X, x_1) +
\left( -Q(X,x_1) + (X-x_1) \frac{Q(X, x_1)-Q(X, x_0)}{x_0-x_1} \right) \, w_0.
\end{multline*}
Then the conclusion follows from the boundedness of $Q$ and $Q'$ and the mean value theorem.
\end{proof*}

\section{Comparison with previous results}
\label{compprev}

In \cite{TraoTh}, we had used a different definition for a Lempert function with multiplicities. We state it with the same notations as in Definition \ref{defLempert}.
\begin{defn}
\label{prevdef}
Given a system $S$ as
in Definition \ref{defgreen},
we write $\tau_j := \tau_{\Psi_j}$. 

Let $\varphi \in
Hol(\D,\Omega)$
and $\alpha_j \in \D$, $1\le j \le N$.
We say that $(\varphi, (\alpha_j)_{1\le j \le N})$ is \emph{admissible (for $S$, $z$)
in the old sense} if 
\begin{multline*}
\varphi (0)=z, \mbox{ and there exists }
 U_j
\mbox{\rm\ a neighborhood of }\zeta_j \\
\mbox{ s.t. }
\Psi_j (\varphi(\zeta) -a_j) \le \tau_j \log|\zeta-\zeta_j| + C_j,
\forall
\zeta \in U_j, 1\le j \le N.
\end{multline*}
In this case, we write (with the convention that $0 \cdot \infty =0$)
$$
\mathcal S (\varphi, (\alpha_j)_{1\le j \le N}):= \sum_{j=1}^N  \tau_j
\log |\alpha_j|   .
$$
 Then
the \emph{old generalized Lempert function} is defined by
\begin{multline*}
L^\Omega_S (z):=
 L_S (z)
\\
:= \inf
\left\lbrace
\mathcal S (\varphi, (\alpha_j)_{1\le j \le N}) : (\varphi, (\alpha_j)_{1\le j \le N})
\mbox{ is admissible for }S, z \mbox{ in the old sense }
\right\rbrace .
\end{multline*}
\end{defn}

Recall also that since the functional $L$ did not enjoy monotonicity properties, another definition was given in \cite{TraoTh}.
\begin{defn}
 \label{modifLemp}
Let $S:= \{(a_j,\Psi_j): 1 \le j \le N\}$ and $S_1:= \{(a_j,\Psi^1_j):
1 \le j \le N\}$ where $a_j \in \Omega$ and $\Psi_j$, $\Psi^1_j$ are local 
indicators. 
We define
$$
\tilde L_S (z) := \inf \{ L_{S^1}(z) : \Psi^1_j \ge \Psi_j + C_j, 1
\le j \le N\} .
$$
\end{defn}
\begin{lemma}
\label{comparison}
If $S=\{(a_j, \Psi_j), 1\le j \le N\}$, where the $\Psi_j$
are elementary local indicators,
then for any $z \in \Omega$, $\mathcal L_S (z) \le \tilde L_S (z)$.
\end{lemma}
\begin{proof}
Since the functional $\mathcal L$ is monotonic by Theorem
\ref{monotone}, it will be enough to show that  $\mathcal L_S (z) \le  L_S (z)$
for any system $S$. If we have a map $\varphi$ which is admissible in the sense of  Definition \ref{prevdef}, we can take $A_j:=\{\alpha_j\}$, and 
$\Psi_j (\varphi(\zeta) -a_j) \le \tau_j \log|\zeta-\alpha_j| + C_j$
 implies that $m_{\varphi, a_j, \Psi_j} (\alpha_j) \ge \tau_j$, which by
Definition \ref{multphia} means that $m_{\varphi, a_j, \Psi_j} (\alpha_j) = \tau_j$. So that any such $\varphi$ is admissible in the sense of Definition
\ref{defLempert}, and
$$
\mathcal S (\varphi, (a_j)_{1\le j \le N})  = \mathcal S (\varphi, (A_j)_{1\le j \le N}),
$$
and the desired inequality follows.
\end{proof}

We now return to the study of the example presented in \cite{TraoTh}. 
Let us recall the notations.
For $z \in \mathbb D^2$,
$$
\Psi_0 (z) := \max(\log|z_1|,\log|z_2|),
\quad \Psi_V (z) := \max(\log|z_1|,2\log|z_2|).
$$
Here $V$ stands for "vertical", for the
obvious
reasons : for $a \in \mathbb D^2$,
\newline
$\Psi_j ( \varphi (\zeta) -a) \le \tau_j \log|\zeta -\zeta_0| + C$
translates to ($\tau_0=1$, $\tau_V=2$):
\begin{eqnarray*}
\varphi (\zeta_0) = a, &\mbox{ when }& j =0, \\
\varphi (\zeta_0) = a, \varphi'_1(\zeta_0) =0 &\mbox{ when }& j =V .
\end{eqnarray*}
For $a$, $b \in \mathbb D$ and $\eps \in \mathbb C$, let
\begin{eqnarray*}
S_\eps &:=& \{ ((a,0), \Psi_0);((b,0), \Psi_0);((b,\eps), \Psi_0);((a,\eps), \Psi_0) \}
\\
S &:=& \{ ((a,0), \Psi_V);((b,0), \Psi_V) \}.
\end{eqnarray*}
Those are product set situations, and the Green functions are explicitly known. For $w \in \mathbb D$, denote by $\phi_w$ the unique involutive holomorphic
automorphism of the disk which exchanges $0$ and $w$:
$$
\phi_w (\zeta) := \frac{w-\zeta}{1-\zeta \bar w}.
$$
Then
\begin{multline*}
G_S (z_1,z_2) = \max \left( \log |\phi_a(z_1) \phi_b(z_2)|,
2 \log |z_2| \right), \\
G_{S_\eps}  (z_1,z_2) = \max \left( \log |\phi_a(z_1) \phi_b(z_2)|,
 \log |z_2 \phi_\eps(z_2)| \right).
\end{multline*}
The following is proved in \cite[p. 397]{TraoTh}.
\begin{prop}
\label{basicineq}
If $b=-a$ and $|a|^2<|\gamma|<|a|$, then
$G_S (0,\gamma)<\tilde L_S (0,\gamma)$.
\end{prop}
It follows from our Theorem \ref{conv} that for any $z \in \mathbb D^2$,
$\lim_{\eps\to 0} L_{S_\eps} (z) = \mathcal L_S(z)$, and in particular, using Lemma
\ref{comparison}, we find again the result laboriously
obtained in \cite[Proposition 6.1]{TTppt}:
$\limsup_{\eps\to 0} L_{S_\eps} (z) \le \tilde L_S(z)$.
It is a consequence of \cite[Theorem 5.1]{TraoTh} (or equivalently
\cite[Theorem 6.2]{TTppt}) that for $b=-a$ and $|a|^{3/2} < |\gamma| <|a|$,
then $\mathcal L_S (0,\gamma) > G_S (0,\gamma) $;
the motivation then was to obtain the counterexample $ L_{S_\eps} (0,\gamma) >
 G_{S_\eps} (0,\gamma)$ for $|\eps|$ small enough.

On the other hand, when $ |\gamma| < |a|^{3/2}$, the old generalized Lempert function
doesn't provide the correct limit of the single pole Lempert functions. 
\begin{prop}
\label{distinct}
For $b=-a$ and $|a|^2 < |\gamma| < |a|^{3/2}$,
$\mathcal L_S (0,\gamma) < \tilde L_S (0,\gamma) $.
\end{prop}
\begin{proof}
Since Proposition \ref{basicineq} implies that
$\tilde L_S (0,\gamma) > G_S (0,\gamma) = 2 \log |a|$,
it will be enough to provide a mapping $\varphi$
and sets $A_1, A_2$  admissible in the sense of
Definition \ref{defLempert} such that
$\mathcal S (\varphi; A_1, A_2) \le 2 \log |a|$.
We restrict ourselves to $a>0$. We now choose $A_1:=\{\zeta_1, \zeta_4\}$,
$A_2:=\{\zeta_2\}$, with
$$
\zeta_2:= \sqrt{a}, \quad
\zeta_1 := \phi_{\zeta_2} \left( \sqrt{\frac{2a}{1+a^2}} \right), \quad
\zeta_4 := \phi_{\zeta_2} \left( -\sqrt{\frac{2a}{1+a^2}} \right),
$$
and
$$
\varphi_1 (\zeta)
:= \phi_{-a} \left( - \phi_{\zeta_2} (\zeta)^2 \right),
\varphi_2 (\zeta)
:= \frac{\gamma}{\zeta_1\zeta_2\zeta_4} \phi_{\zeta_1} (\zeta) \phi_{\zeta_2} (\zeta) \phi_{\zeta_4} (\zeta).
$$
From those definitions it is clear that $\varphi_1 (\mathbb D) \subset \mathbb D$
and that
$$
\varphi_1 (\zeta_2) = -a,  \varphi_1' (\zeta_2) =0; \quad
\varphi_2 (\zeta_j) = 0, \mbox{ for } j = 1, 2, 4.
$$
 Furthermore, using the involutivity of $\phi_{\zeta_2} $,
$$
\varphi_1 (\zeta_1) = \varphi_1 (\zeta_4) =  \phi_{-a} \left(- \frac{2a}{1+a^2}\right)
=   \phi_{-a} \left( \phi_{-a} (a) \right) = a.
$$
So the map $\varphi$ hits the poles, and
$$
m_{\varphi, (a,0), \Psi_V} (\zeta_1) \ge 1, \quad
m_{\varphi, (a,0), \Psi_V} (\zeta_4) \ge 1, \quad
m_{\varphi, (-a,0), \Psi_V} (\zeta_2) = 2.
$$
To see that actually $m_{\varphi, (a,0), \Psi_V} (\zeta_j)=1$, for $j=1, 4$, 
notice that, since $\varphi_1$ only admits one critical point, $\zeta_2$, 
and since $\zeta_1\neq \zeta_2$ and $\zeta_4\neq \zeta_2$, we
must have $ \varphi_1' (\zeta_j)\neq 0$, $j=1, 4$.

Thus $\varphi$ is admissible in the sense of Definition \ref{defLempert},
and
$$
\mathcal S (\varphi; A_1, A_2) = \log |\zeta_1| + 2 \log |\zeta_2| +\log |\zeta_4| =
\log |\zeta_1 \zeta_4 \zeta_2^2|.
$$
We need to compute
\begin{multline*}
\zeta_1 \zeta_4=
\phi_{\sqrt a} \left( \sqrt{\frac{2a}{1+a^2}} \right) \cdot
\phi_{\sqrt a} \left( -\sqrt{\frac{2a}{1+a^2}} \right) \\
= \frac{\sqrt a -  \sqrt{\frac{2a}{1+a^2}} }{1-\sqrt a \sqrt{\frac{2a}{1+a^2}}}
\cdot
\frac{\sqrt a +  \sqrt{\frac{2a}{1+a^2}} }{1 + \sqrt a \sqrt{\frac{2a}{1+a^2}}}
=\frac{ a -  \frac{2a}{1+a^2} }{1- a   \frac{2a}{1+a^2} }
= \phi_a \left( \phi_a (-a) \right) = -a.
\end{multline*}
From this we deduce $|\zeta_1 \zeta_2 \zeta_4| = a^{3/2} > |\gamma|$, and
therefore $\varphi_2 (\mathbb D) \subset \mathbb D$; and
$ |\zeta_1 \zeta_4 \zeta_2^2| = a^2$, therefore
$$
\mathcal S (\varphi; A_1, A_2) \le \log |\zeta_1 \zeta_4 \zeta_2^2| =
2 \log |a|,
$$
q.e.d.
\end{proof}
\bibliographystyle{amsplain}

\vskip.2cm

Pascal J. Thomas

Institut de Math\'ematiques de Toulouse

CNRS UMR 5219

UFR MIG

Universit\'e Paul Sabatier

F-31062 TOULOUSE CEDEX 9

France

pthomas@math.univ-toulouse.fr

\vskip.5cm

Nguyen Van Trao

Department of Mathematics

Dai Hoc Su Pham 1 (Pedagogical Institute of Hanoi)

Cau Giay, Tu Liem

Ha Noi

Viet Nam

ngvtrao@yahoo.com

\end{document}